\documentclass[a4paper,12pt]{amsart} 

\usepackage[ansinew]{inputenc} 
\usepackage[T1]{fontenc}
\usepackage[english]{babel} 

\usepackage{enumerate} 
\usepackage{booktabs}

\usepackage{amsmath} 
\usepackage{amsfonts} 
\usepackage{amssymb}
\usepackage{amsthm}
\usepackage{latexsym}
\usepackage{bbm}
\usepackage[top=3cm , bottom=3cm , left=3cm , right=3cm]{geometry}
\usepackage{esint}

\makeatletter

\@namedef{subjclassname@2010}{%

  \textup{2010} Mathematics Subject Classification}

\makeatother

\theoremstyle{plain}
\newtheorem{lem}{Lemma} 
\newtheorem{thm}[lem]{Theorem}

\theoremstyle{definition} 
 
\newtheorem*{definition*}{Definition}

\theoremstyle{remark}
 
\newcommand{\R}{\mathbbm{R}}

\newcommand{\Z}{\mathbbm{Z}}

\newcommand{\E}{\mathbbm{E}}
\newcommand{\SP}{\mathbbm{S}}
\newcommand{\prob}{\mathbbm{P}}
\newcommand{\D}{\,\textup{d}}
\newcommand{\la}{\langle}
\newcommand{\ra}{\rangle}
\newcommand{\supp}{\text{supp}\,}

\title[Vector-valued tent spaces $T^1$ and $T^{\infty}$]{The vector-valued tent spaces $T^1$ and $T^{\infty}$}
\author[M. Kemppainen]{Mikko Kemppainen}
\address{Department of Mathematics and Statistics, University of Helsinki,
Gustaf Hällströmin katu 2b, FI-00014 Helsinki, Finland}

\email{mikko.k.kemppainen@helsinki.fi}

\begin{document}




\maketitle

\begin{abstract}
  Tent spaces of vector-valued functions were recently studied by
  Hytönen, van Neerven and Portal with an eye on applications to
  $H^{\infty}$-functional calculi. This paper extends their results to the endpoint
  cases $p=1$ and $p=\infty$ along the lines of earlier work by Harboure, Torrea and Viviani
  in the scalar-valued case.
  The main result of the paper is an atomic decomposition in the case $p=1$, which relies on a new geometric
  argument for cones. A result on the duality of these spaces is also given.
\end{abstract}

\section{Introduction}\label{s:1}

Coifman, Meyer and Stein introduced in \cite{CMSTENTSPACES} the concept of tent spaces that provides a
neat framework for several ideas and techniques in Harmonic Analysis. In particular, they defined
the spaces $T^p$, $1 \leq p < \infty$, 
that are relevant for square functions, and consist of functions $f$ on 
the upper half-space $\R^{n+1}_+$ for which the $L^p$ norm
of the conical square function is finite:
\begin{equation*}
  \int_{\R^n} \Big( \int_{\Gamma (x)} |f(y,t)|^2 \frac{\D y \D t}{t^{n+1}} \Big)^{p/2} \D x < \infty ,
\end{equation*}
where $\Gamma (x)$ denotes the cone $\{ (y,t)\in\R^{n+1}_+ : | x - y | < t \}$ at $x\in\R^n$.
Typical functions in these spaces arise for instance from harmonic extensions $u$ to
$\R^{n+1}_+$ of $L^p$ functions on $\R^n$
according to the formula $f(y,t) = t\partial_tu (y,t)$.

Tent spaces were approached by Harboure, Torrea and Viviani in \cite{HARBOURE} as $L^p$ spaces of
$L^2$-valued functions, which gave an abstract way to deduce many of their basic properties. Indeed,
for $1 < p < \infty$, the mapping $Jf(x) = 1_{\Gamma (x)}f$
is readily seen to embed $T^p$ in $L^p(\R^n ; L^2(\R^{n+1}_+))$, when
$\R^{n+1}_+$ is equipped with the measure $\D y \D t / t^{n+1}$. Furthermore, they showed that
$T^p$ is embedded as a complemented subspace, which not only implies its completeness, but also
gives a way to prove a few other properties, such as equivalence of norms defined by cones of different
aperture and the duality $(T^p)^* \simeq T^{p'}$, where $1/p + 1/p' = 1$.

Treatment of the endpoint cases $p=1$ and $p=\infty$ requires more careful inspection. Firstly,
the space $T^{\infty}$ was defined in \cite{CMSTENTSPACES} as the space of
functions $g$ on $\R^{n+1}_+$ for which
\begin{equation*}
  \sup_B \frac{1}{|B|} \int_{\widehat{B}} |g(y,t)|^2 \frac{\D y \D t}{t} < \infty ,
\end{equation*}
where the supremum is taken over all balls $B\subset\R^n$ and where 
$\widehat{B} \subset \R^{n+1}_+$ denotes the ``tent''
over $B$ (see Section \ref{s:4}). 
The tent space duality is now extended
to the endpoint case as $(T^1)^* \simeq T^{\infty}$. Moreover, functions in $T^1$ admit a decomposition into
atoms $a$ each of which is supported in $\widehat{B}$ for some ball $B\subset \R^n$ and satisfies
\begin{equation*}
  \int_{\widehat{B}} |a(y,t)|^2 \frac{\D y \D t}{t} \leq \frac{1}{|B|} .
\end{equation*}
As for the embeddings, it is proven in \cite{HARBOURE} that $T^1$ embeds in the 
$L^2(\R^{n+1}_+)$-valued Hardy space
$H^1(\R^n ; L^2(\R^{n+1}_+))$, while $T^{\infty}$ embeds in $\textup{BMO} (\R^n ; L^2(\R^{n+1}_+))$ --
the space of $L^2(\R^{n+1}_+)$-valued functions with bounded mean oscillation.

The study of vector-valued analogues of these spaces was initiated by Hyt\"onen, van Neerven and Portal in
\cite{HVNPCONICAL}, where they followed the ideas from \cite{HARBOURE} and proved the
analogous embedding results for $T^p(X)$ with $1 < p < \infty$ under the assumption that $X$ is UMD. 
It should be noted that, for non-Hilbertian
$X$, the $L^2$ integrals had to be replaced by stochastic integrals or some equivalent
objects, which in turn required further adjustments in proofs, namely the lattice maximal functions
that appeared in \cite{HARBOURE} were replaced by an appeal to Stein's inequality for conditional
expectation operators. Later on, Hyt\"onen and Weis 
provided in \cite{HYTONENWEISPARAPRODUCTS} a scale of vector-valued
versions of the quantity appearing above in the definition of $T^\infty$. 

This paper continues the work on the endpoint cases and provides definitions for $T^1(X)$
and $T^{\infty}(X)$. The main result decomposes a $T^1(X)$ function into atoms
using a geometric argument for cones. The original decomposition argument in \cite{CMSTENTSPACES} is
inherently scalar-valued and not as such suitable for stochastic integrals. Moreover, the spaces
$T^1(X)$ and $T^{\infty}(X)$ are embedded in certain Hardy and \textup{BMO} spaces, respectively, 
much in the spirit of \cite{HARBOURE}. 
The theory of vector-valued stochastic integration (see van Neerven and Weis \cite{JVNSTOCHINT})
is used throughout the paper.

\subsection*{Acknowledgements}

I gratefully acknowledge the support from
the Finnish National Graduate School in Mathematics and its Applications and from the Academy of Finland,
grant 133264.
I would also like to thank Tuomas Hyt\"onen, Jan van Neerven, Hans-Olav Tylli and Mark Veraar for insightful comments and
conversations.

\section{Preliminaries}\label{s:2}

\subsection*{Notation.}

Random variables are taken to be defined on a fixed probability space whose 
probability measure and expectation are denoted by 
$\prob$ and $\E$. The integral average
(with respect to Lebesgue measure) over a measurable set $A\subset\R^n$ is 
written as $\fint_A = |A|^{-1}\int_A$, where
$|A|$ stands for the Lebesgue measure of $A$. For a ball $B$ in $\R^n$ we write $x_B$ and $r_B$ for its
center and radius, respectively.
Throughout the paper $X$ is 
assumed to be a real Banach space and 
$\la \xi , \xi^* \ra$ is used to denote the duality pairing between $\xi\in X$ and $\xi^*\in X^*$.
Isomorphism of Banach spaces is expressed using $\simeq$.
By $\alpha \lesssim \beta$ it is meant that there exists a constant $C$ such that
$\alpha \leq C\beta$. Quantities $\alpha$ and $\beta$ are comparable, $\alpha \eqsim \beta$, if
$\alpha \lesssim \beta$ and $\beta \lesssim \alpha$.

\subsection*{Stochastic integration.}

We start by discussing the correspondence between
Gaussian random measures and stochastic integrals of real-valued functions.
Recall that a Gaussian random measure on a $\sigma$-finite measure space
$(M,\mu)$ is a mapping $W$ that takes
subsets of $M$ with finite measure to (centered) Gaussian random variables
in such a manner that 
\begin{itemize}
\item the variance $\E W(A)^2 = \mu (A)$,
\item for all disjoint $A$ and $B$
the random variables $W(A)$ and $W(B)$ are independent and
$W(A\cup B) = W(A) + W(B)$ almost surely.
\end{itemize}
Since for Gaussian random variables the notions of
independence and orthogonality are equivalent, it suffices to consider their pairwise independence
in the definition above.
Given a Gaussian random measure $W$, we obtain a linear isometry
from $L^2(M)$ to $L^2(\prob )$ 
-- our stochastic integral --
by first defining $\int_M 1_A \D W = W(A)$ and then extending by linearity and density to the whole
of $L^2(M)$. On the other hand, if we are in possession of such an isometry, we may define
a Gaussian random measure $W$ by sending a subset $A$ of $M$ with finite measure to the stochastic
integral of $1_A$. For more details, see Janson \cite[Chapter 7]{JANSON}.

A function $f: M \to X$ is said to be weakly $L^2$ if $\la f(\cdot ) , \xi^* \ra$ is in $L^2(M)$ for all
$\xi^*\in X^*$. Such a function is said to be 
\emph{stochastically integrable} (with respect to a Gaussian random
measure $W$) if there exists a (unique) random variable $\int_M f \D W$ in $X$ so that for all $\xi^*\in X^*$ we have
\begin{equation*}
  \Big\la \int_M f \D W , \xi^* \Big\ra = \int_M \la f(t) , \xi^* \ra \D W(t) \quad \text{almost surely}.
\end{equation*}
We also say that a function $f$ is stochastically
integrable over a measurable subset $A$ of $M$ if $1_Af$ is stochastically integrable.
Note in particular that each function $f = \sum_k f_k \otimes \xi_k$ in 
the algebraic tensor product $L^2(M) \otimes X$ is stochastically integrable and that
\begin{equation*}
  \int_M f \D W = \sum_k \Big( \int_M f_k \D W \Big) \xi_k .
\end{equation*}
A detailed theory of vector-valued stochastic integration can be found in 
van Neerven and Weis \cite{JVNSTOCHINT}, see also Rosi\'nski and Suchanecki \cite{ROSINSKI}.
Stochastic integrals have a number of nice properties (see \cite{JVNSTOCHINT}):
\begin{itemize}
\item Khintchine--Kahane inequality: For every stochastically integrable $f$ we have
      \begin{equation*}
        \Big( \E \Big\| \int_M f \D W \Big\|^p \Big)^{1/p} 
        \eqsim \Big( \E \Big\| \int_M f \D W \Big\|^q \Big)^{1/q}
      \end{equation*}
      whenever $1\leq p,q < \infty$.
\item Covariance domination: 
      If a function $g\in L^2(M) \otimes X$ is dominated by a function
      $f\in L^2(M) \otimes X$ in covariance, that is, if
      \begin{equation*}
        \int_M \la g(t),\xi^* \ra ^2 \D\mu (t) \leq \int_M \la f(t),\xi^* \ra ^2 \D\mu (t)
      \end{equation*}
      for all $\xi^*\in X^*$, then
      \begin{equation*}
        \E \Big\| \int_M g \D W \Big\|^2 \leq \E \Big\| \int_M f \D W \Big\|^2 .
      \end{equation*} 
\item Dominated convergence: 
      If a sequence $(f_k)$ of stochastically integrable functions is dominated in covariance
      by a single stochastically integrable function and
      \begin{equation*}
        \int_M \la f_k(t),\xi^* \ra ^2 \D\mu (t) \to 0
      \end{equation*}
      for all $\xi^*\in X^*$, then 
      \begin{equation*}      
        \E \Big\| \int_M f_k \D W \Big\|^2 \to 0 .
      \end{equation*}
      In particular,
      if a sequence $(A_k)$ of measurable sets satisfies $1_{A_k} \to 0$ pointwise almost everywhere, then
      for every $f$ in $L^2(M) \otimes X$ we have
      \begin{equation*}
        \E \Big\| \int_{A_k} f \D W \Big\|^2 \to 0 .
      \end{equation*}
\end{itemize}

The expression
\begin{equation*}
   \Big( \E \Big\| \int_M f \D W \Big\|^2 \Big)^{1/2}
\end{equation*}
defines a norm on the space of (equivalence classes of) strongly measurable stochastically integrable functions
$f:M\to X$. However, the norm is not generally complete, unless $X$ is a Hilbert space. For convenience, we operate mainly
with functions in $L^2(M) \otimes X$ and denote their completion under the norm above by $\gamma (M;X)$.

This space can be identified with the space of $\gamma$-radonifying operators from
$L^2(M)$ to $X$ (see \cite{JVNSTOCHINT} and the survey \cite{GAMMARAD}). 
We note the following facts:
\begin{itemize}
\item Given an $m\in L^{\infty}(M)$, the multiplication operator $f\mapsto mf$ on $L^2(M) \otimes X$ has norm
$\| m \|_{L^{\infty}(M)}$.
\item For K-convex $X$ (see \cite[Section 10]{GAMMARAD}) the duality $\gamma (M;X)^* = \gamma (M;X^*)$
holds and realizes for $f\in L^2(M) \otimes X$ and $g\in L^2(M) \otimes X^*$ via
\begin{equation*}
  \la f , g \ra = \int_M \la f(t),g(t) \ra \D\mu (t) .
\end{equation*}
\end{itemize}

A family $\mathcal{T}$ of operators in $\mathcal{L}(X)$ is said to be \emph{$\gamma$-bounded} if
for every finite collection of operators $T_k\in\mathcal{T}$ and vectors $\xi_k\in X$ we have
\begin{equation*}
  \E \Big\| \sum_k \gamma_k T_k \xi_k \Big\|^2 \lesssim \E \Big\| \sum_k \gamma_k \xi_k \Big\|^2 ,
\end{equation*}
where $(\gamma_k)$ is an independent sequence of standard Gaussians.

Observe, that families of operators obtained by composing operators from (a finite number of)
other $\gamma$-bounded
families are also $\gamma$-bounded. It follows from covariance domination and Fubini's theorem, 
that the family of
operators $f\mapsto mf$ is $\gamma$-bounded on $L^p(\R^n ;X)$ whenever the multipliers $m$ are chosen from
a bounded set in $L^{\infty}(\R^n )$.

The following continuous-time result for $\gamma$-bounded families is folklore (to be found in Kalton and 
Weis \cite{WEISUNPUB}):

\begin{lem}
\label{gammabounded}
Assume that $X$ does not contain a closed subspace isomorphic to $c_0$.
If the range of an $X$-strongly measurable function $A:M\to \mathcal{L}(X)$ 
is $\gamma$-bounded, then for every strongly measurable stochastically integrable function $f: M\to X$ the 
strongly measurable function $t\mapsto A(t)f(t) : M \to X$ is also stochastically integrable and satisfies
\begin{equation*}
  \E \Big\| \int_M A(t)f(t) \D W(t) \Big\| ^2 \lesssim \E \Big\| \int_M f(t) \D W(t) \Big\|^2 .
\end{equation*}
\end{lem}

Recall that $X$-strong measurability of a function $A:M\to \mathcal{L}(X)$ requires
$A(\cdot )\xi : M \to X$ to be strongly measurable for every $\xi\in X$.
For simple functions $A:M\to\mathcal{L}(X)$ the lemma above is immediate from the definition of $\gamma$-boundedness
and requires no assumption regarding containment of $c_0$,
as the function $t\mapsto A(t)f(t) : M \to X$ is also in $L^2(M)\otimes X$. Assuming $A$ to be simple is anyhow too
restrictive for applications and to consider nonsimple functions $A$ 
we need to handle more general stochastically 
integrable functions than just those in $L^2(M) \otimes X$.

Our choice of $(M,\mu)$ will be the upper half-space
$\R^{n+1}_+ = \R^n \times (0,\infty )$ equipped with the measure $\D y \D t / t^{n+1}$.
We will simplify our notation and write $\gamma (X) = \gamma (\R^{n+1}_+;X)$
-- in what follows, stochastic integration is performed on $\R^{n+1}_+$.

\subsection*{The UMD-property and averaging operators.}

It is often necessary to assume that our Banach space $X$ is UMD.
This has the crucial 
implication, known as \emph{Stein's inequality} (see Bourgain \cite{BOURGAINSTEININEQ} and 
Cl\'ement et al. \cite{CLEMENT}), 
that every increasing family of conditional expectation operators is $\gamma$-bounded on
$L^p(X)$ whenever $1 < p < \infty$. 
While this is proven in the given references only in the case of probability spaces,
it can be generalized to the $\sigma$-finite case such as ours with no difficulty.
Namely, let us consider filtrations on $\R^n$ generated by
systems of dyadic cubes, that is, by collections $\mathcal{D} = \bigcup_{k\in\Z} \mathcal{D}_k$, where
each $\mathcal{D}_k$ is a disjoint cover of $\R^n$ consisting of cubes $Q$ of the form $x_Q + [0,2^{-k})^n$
and every $Q\in\mathcal{D}_k$ is a union of $2^n$ cubes in $\mathcal{D}_{k+1}$.
The conditional expectation operators or averaging operators are then given for each integer $k$ by
\begin{equation*}
  f \mapsto \sum_{Q\in\mathcal{D}_k} 1_Q \fint_Q f , \quad f\in L^1_{\text{loc}}(\R^n ; X) .
\end{equation*}
Composing such an operator with multiplication by an indicator $1_Q$ of a dyadic cube $Q$, we arrive
through Stein's inequality to the conclusion that the family $\{ A_Q \}_{Q\in\mathcal{D}}$ of 
localized averaging operators
\begin{equation*}
  A_Q f = 1_Q \fint_Q f ,
\end{equation*}
is $\gamma$-bounded on $L^p(\R^n ; X)$ whenever $1 < p < \infty$.
The following result of Mei \cite{TAOMEI} allows us to replace dyadic cubes by balls:

\begin{lem}
  There exist $n+1$ systems of dyadic cubes such that every ball
  $B$ is contained in a dyadic cube $Q_B$ from one of the systems and $|B| \lesssim |Q_B|$.
\end{lem}
Stein's inequality together with the lemma above guarantees that the family 
$\{ A_B : B \text{ ball in } \R^n \}$ is $\gamma$-bounded on $L^p(\R^n ; X)$ whenever $1 < p < \infty$. 
Indeed, for each ball $B$ we can write
\begin{equation*}
  A_B = 1_B \frac{|Q_B|}{|B|} A_{Q_B} 1_B .
\end{equation*}
This was proven already in \cite{HVNPCONICAL}.

It will be useful to consider smoothed or otherwise different versions of indicators
$1_B(x) = 1_{[0,1)} ( |x - x_B| / r_B )$. Given a measurable $\psi : [0,\infty ) \to \R$
with $1_{[0,1)} \leq |\psi | \leq 1_{[0,\alpha )}$ for some $\alpha > 1$, 
we define the averaging operators
\begin{equation*}
  A_{y,t}^{\psi}f(x) = \psi \Big( \frac{|x-y|}{t} \Big) \frac{1}{c_{\psi} t^n} 
  \int_{\R^n} \psi \Big( \frac{|z-y|}{t} \Big) f(z) \D z ,
  \quad x\in\R^n ,
\end{equation*}
where
\begin{equation*}
  c_\psi = \int_{\R^n} \psi (|x|)^2 \D x .
\end{equation*}
Again, under the assumption that $X$ is UMD and $1 < p < \infty$,
the $\gamma$-boundedness of the 
family $\{ A_{y,t}^{\psi} : (y,t)\in\R^{n+1}_+ \}$
of operators on $L^p(\R^n ; X)$ follows at once when we write
\begin{equation*}
  A_{y,t}^{\psi} = \psi \Big( \frac{|\cdot - y|}{t} \Big) 
  \frac{|Q_{B(y,\alpha t)}|}{c_{\psi} t^n} A_{Q_{B(y,\alpha t)}}
  \psi \Big( \frac{|\cdot - y|}{t} \Big) .
\end{equation*}

Observe, that the 
function $(y,t)\mapsto A_{y,t}^\psi$ from
$\R^{n+1}_+$ to $\mathcal{L}(L^p(\R^n ; X))$
is $L^p(\R^n ; X)$-strongly measurable.
Recall also that every UMD space is K-convex and cannot contain a closed subspace isomorphic to $c_0$.

\section{Overview of tent spaces}\label{s:3}

\subsection*{Tent spaces $T^p(X)$.}

Let us equip the upper half-space $\R^{n+1}_+$ with the measure $\D y \D t / t^{n+1}$
and a Gaussian random measure $W$.
Recall the definition of the cone $\Gamma (x) = \{ (y,t)\in\R^{n+1}_+ : | x - y | < t \}$ at $x\in\R^n$.

Let $1\leq p < \infty$.
We wish to define a norm on the space of 
functions $f:\R^{n+1}_+ \to X$ for which $1_{\Gamma (x)}f \in L^2(\R^{n+1}_+) \otimes X$ for almost every $x\in\R^n$ by
\begin{equation*}
  \| f \|_{T^p(X)} 
  = \Big( \int_{\R^n} \Big( \E \Big\| 
  \int_{\Gamma (x)} f \D W \Big\| ^2 \Big) ^{p/2} \D x \Big)^{1/p}
\end{equation*}
and use the Khintchine--Kahane inequality to write
\begin{equation*}
  \| f \|_{T^p(X)} 
  \eqsim \Big( \E \Big\| \int_{\Gamma (\cdot )} f \D W \Big\|_{L^p(\R^n ; X)}^p \Big)^{1/p} ,
\end{equation*}
but issues concerning measurability need closer inspection. 

\begin{lem}
\label{measurability}
  Suppose that $f:\R^{n+1}_+ \to X$ is such that $1_{\Gamma (x)}f \in L^2(\R^{n+1}_+) \otimes X$ 
  for almost every $x\in\R^n$. Then 
  \begin{enumerate}  
  \item the function $x\mapsto 1_{\Gamma(x)} f$ is strongly measurable from $\R^n$ to $\gamma (X)$.
  \item the function $x \mapsto \int_{\Gamma (x)} f \D W$ is strongly measurable
        from $\R^n$ to $L^2(\prob ; X)$ and may be considered, when 
        $\| f \|_{T^p(X)} < \infty$, as a random $L^p(\R^n ; X)$ function.
  \item the function $x\mapsto ( \E \| \int_{\Gamma (x)} f \D W \|^2 )^{1/2}$ agrees almost everywhere
        with a lower semicontinuous function so that the set
        \begin{equation*}
          \Big\{ x \in \R^n : \Big( \E \Big\| \int_{\Gamma (x)} f \D W \Big\|^2 \Big)^{1/2} > \lambda \Big\}
        \end{equation*}
        is open whenever $\lambda > 0$. 
  \end{enumerate}
  \begin{proof}
    Denote by $A_k$ the set $\{ (y,t)\in \R^{n+1}_+ : t > 1/k \}$ 
    and write $f_k = 1_{A_k} f$. It is clear that for 
    each positive integer $k$, the functions
    $x\mapsto 1_{\Gamma (x)}f_k$ and $x\mapsto \int_{\Gamma (x)} f_k \D W$ are 
    strongly measurable and continuous since
    \begin{equation*}
      \E \Big\| \int_{\Gamma (x) \Delta \Gamma (x')} f_k \D W \Big\|^2 \to 0, 
      \quad \text{as} \quad x \to x' .
    \end{equation*}
    Furthermore, $1_{\Gamma (x)}f_k \to 1_{\Gamma (x)}f$ in $\gamma (X)$ for almost every $x\in\R^n$ since
    \begin{equation*}
      \E \Big\| \int_{\Gamma (x)} (f - f_k) \D W \Big\|^2 
      = \E \Big\| \int_{\Gamma(x)\setminus A_k} f \D W \Big\|^2 \to 0 .
    \end{equation*} 
    Consequently, $x\mapsto 1_{\Gamma (x)} f$ and $x\mapsto \int_{\Gamma (x)} f \D W$ 
    are strongly measurable.
    Moreover, the pointwise limit of an increasing sequence of real-valued 
    continuous functions is lower semicontinuous,
    which proves the third claim.
  \end{proof}
\end{lem}

\begin{definition*}
  Let $1 \leq p < \infty$. The tent space $T^p(X)$ is defined as the completion
  under $\| \cdot \|_{T^p(X)}$ of the space of (equivalence classes of)
  functions $\R^{n+1}_+ \to X$ (in what follows, ``$T^p(X)$ functions'') such that 
  $1_{\Gamma (x)}f \in L^2(\R^{n+1}_+) \otimes X$ for almost every $x$ in $\R^n$ and 
  $\| f \|_{T^p(X)} < \infty$.
\end{definition*}

As was mentioned in the previous section, it is useful to consider the more general situation where
the indicator of a ball is replaced by
a measurable function $\phi : [0,\infty ) \to \R$ with
$1_{[0,1)} \leq |\phi | \leq 1_{[0,\alpha )}$ for some $\alpha > 1$. 
Let us assume in addition, that $\phi$ is continuous at $0$.
For functions $f:\R^{n+1}_+ \to X$ such that 
\begin{equation*}
  (y,t)\mapsto \phi (|x-y|/t) f(y,t) \in L^2(\R^{n+1}_+) \otimes X
\end{equation*}  
for almost
every $x\in\R^n$, the strong measurability of
\begin{equation*}
  x \mapsto \Big( (y,t) \mapsto \phi \Big( \frac{|x-y|}{t} \Big) f(y,t) \Big) \quad \text{and} \quad
  x \mapsto \int_{\Gamma (x)} \phi \Big( \frac{|x-y|}{t} \Big) f(y,t) \D W(y,t)
\end{equation*}
are treated as in the case of $\phi (|x-y|/t) = 1_{[0,1)}(|x-y|/t) = 1_{\Gamma (x)}(y,t)$.

\subsection*{Embedding $T^p(X)$ into $L^p(\R^n ; \gamma (X))$.}

A collection of results from the paper \cite{HVNPCONICAL} by Hyt\"onen, van Neerven and Portal 
is presented next. 
Following the idea of Harboure, Torrea and Viviani \cite{HARBOURE},
the tent spaces are embedded into $L^p$ spaces of $\gamma (X)$-valued functions by
\begin{equation*}
  Jf(x) = 1_{\Gamma (x)}f, \quad x\in\R^n .
\end{equation*}
Furthermore, for simple $L^2(\R^{n+1}_+) \otimes X$ -valued functions $F$ on $\R^n$ 
we define an operator $N$ by
\begin{equation*}
  (N F)(x;y,t) = 1_{B(y,t)}(x) \fint_{B(y,t)} F(z;y,t) \D z, \quad x\in\R^n, (y,t)\in\R^{n+1}_+ .
\end{equation*}
Assuming that $X$ is UMD, we can now view $T^p(X)$ as a complemented subspace of $L^p(\R^n ; \gamma (X))$:

\begin{thm}
\label{projection}
  Suppose that $X$ is UMD and let $1 < p < \infty$. Then $N$ extends to a bounded projection on $L^p(\R^n ; \gamma (X))$
  and $J$ extends to an isometry from $T^p(X)$ onto the image of $L^p(\R^n ; \gamma (X))$ under $N$.
\end{thm}

The following result shows the comparability of different tent space norms: 

\begin{thm}
\label{normequiv}
  Suppose that $X$ is UMD, let $1 < p < \infty$ and let $1_{[0,1)} \leq |\phi | \leq 1_{[0, \alpha )}$. 
  For every function $f$ in $T^p(X)$ the function $(y,t)\mapsto \phi (|x-y|/t)f(y,t)$ is stochastically
  integrable for almost every $x\in\R^n$ and  
  \begin{equation*}
    \int_{\R^n} \E \Big\| \int_{\R^{n+1}_+} \phi \Big( \frac{|x-y|}{t} \Big) f(y,t) \D W(y,t) \Big\|^p \D x
    \eqsim \int_{\R^n} \E \Big\| \int_{\Gamma (x)} f \D W \Big\|^p \D x .
  \end{equation*}
\end{thm}

The proof relies on the boundedness of modified projection operators
\begin{equation*}
  (N_\phi F)(x;y,t) 
  = \phi \Big( \frac{|x-y|}{t} \Big) \fint_{B(y,t)} F(z;y,t) \D z, \quad x\in\R^n, (y,t)\in\R^{n+1}_+ . 
\end{equation*}
and the observation that the embedding
\begin{equation*}
  J_\phi f(x;y,t) = \phi \Big( \frac{|x-y|}{t} \Big) f(y,t), \quad x\in\R^n, (y,t)\in\R^{n+1}_+ .
\end{equation*}
can be written as $J_\phi f = N_\phi Jf$.
In particular, this shows that the norms given by cones of different apertures are comparable.
Indeed, choosing $\phi = 1_{[0,\alpha )}$ gives the norm where $\Gamma (x)$ is
replaced by the cone $\Gamma_\alpha (x) = \{ (y,t)\in\R^{n+1}_+ : |x-y| < \alpha t \}$ with aperture
$\alpha > 1$. 

Indentification of tent spaces $T^p(X)$ with complemented subspaces of $L^p(\R^n ; \gamma (X))$ gives a powerful way to deduce
their duality:

\begin{thm}
\label{Tpduality}
  Suppose that $X$ is UMD and let $1 < p < \infty$. 
  Then the dual of $T^p(X)$ is $T^{p'}(X^*)$, where $1/p + 1/p' = 1$, and the duality is realized
  for functions $f\in T^p(X)$ and $g\in T^{p'}(X^*)$ via
  \begin{equation*}
    \la f,g \ra = c_n \int_{\R^{n+1}_+} \la f(y,t),g(y,t) \ra \frac{\D y \D t}{t} ,
  \end{equation*}
  where $c_n$ is the volume of the unit ball in $\R^n$.
\end{thm}

The following theorem combines results from \cite[Theorem 4.8]{HVNPCONICAL} and 
\cite[Corollary 4.3, Theorem 1.3]{HYTONENWEISPARAPRODUCTS}. 
The tent space $T^\infty (X)$ is defined in the next section.

\begin{thm}
  Suppose that $X$ is UMD and let $\Psi$ be a Schwartz function with vanishing integral. Then the operator
  \begin{equation*}
    T_\Psi f(y,t) = \Psi _t \ast f(y)
  \end{equation*}
  is bounded from $L^p(\R^n ; X)$ to $T^p(X)$ whenever $1 < p < \infty$, from 
  $H^1(\R^n ; X)$ to $T^1(X)$ and from
  $\textup{BMO}(\R^n ; X)$ to $T^{\infty}(X)$.
\end{thm}

\section{Tent spaces $T^1(X)$ and $T^{\infty} (X)$}\label{s:4}

Having completed our overview of tent spaces $T^p(X)$ with $1 < p < \infty$ we turn to the
endpoint cases $p=1$ and $p=\infty$, of which the latter remains to be defined. 
As for the case $p=1$,
our aim is to show that $T^1(X)$ is isomorphic to a complemented subspace of 
the Hardy space $H^1(\R^n ; \gamma (X))$ of $\gamma (X)$-valued functions on $\R^n$.
In the case $p=\infty$, we introduce the space 
$T^{\infty}(X)$, which is shown to
embed in $\textup{BMO} (\R^n ; \gamma (X))$, that is, the space of $\gamma (X)$-valued
functions whose mean oscillation is bounded.
The idea of these embeddings was originally put forward by Harboure et al. in the scalar-valued case
(see \cite{HARBOURE}). 

Recall that the tent over an open set $E\subset \R^n$ is defined by
$\widehat{E} = \{ (y,t) \in \R^{n+1}_+ : B(y,t) \subset E \}$ or equivalently by
\begin{equation*}
  \widehat{E} = \R^{n+1}_+ \setminus \bigcup_{x\not\in E} \Gamma (x) .
\end{equation*}
Observe that while cones are open, tents are closed. Truncated cones are also needed:
For $x\in\R^n$ and $r > 0$ we define $\Gamma (x;r) = \{ (y,t)\in\Gamma (x) : t < r \}$.

In \cite{HYTONENWEISPARAPRODUCTS} Hyt\"onen and Weis adjusted the quantities that define
scalar-valued atoms and $T^{\infty}$ functions in terms of tents to more suitable ones that rely on
averages of square functions. More precisely
for scalar-valued $g$ on $\R^{n+1}_+$ we have
\begin{align*}
  \int_B \int_{\Gamma (x;r_B)} |g(y,t)|^2 \frac{\D y \D t}{t^{n+1}} \D x
  &= \int_B \int_{\R^n \times (0,r_B)} 1_{B(y,t)}(x) |g(y,t)|^2 \frac{\D y \D t}{t^{n+1}} \D x \\
  &= \int_0^{r_B} \int_{2B} | g(y,t) |^2 | B \cap B(y,t)|\frac{\D y \D t}{t^{n+1}} ,
\end{align*}
from which one reads
\begin{equation*}
  \int_{\widehat{B}} |g(y,t)|^2 \frac{\D y \D t}{t}
  \lesssim \int_B \int_{\Gamma (x;r_B)} |g(y,t)|^2 \frac{\D y \D t}{t^{n+1}} \D x
  \lesssim \int_{\widehat{3B}} |g(y,t)|^2 \frac{\D y \D t}{t} .
\end{equation*}
This motivates the definition of a \emph{$T^1(X)$ atom} as 
a function $a:\R^{n+1}_+ \to X$ such that for some ball $B$ we have $\supp a \subset \widehat{B}$, 
$1_{\Gamma (x)} a \in L^2(\R^{n+1}_+) \otimes X$ for almost every $x\in B$ and
\begin{equation*}
  \int_B \E \Big\| \int_{\Gamma (x)} a \D W \Big\|^2 \D x \leq \frac{1}{|B|} .
\end{equation*}
Then $1_{\Gamma (x)} a$ differs from zero only when $x\in B$ and so
\begin{equation*}
  \| a \|_{T^1(X)} = \int_{\R^n} \Big( \E \Big\| \int_{\Gamma (x)} a \D W \Big\|^2 \Big)^{1/2} \D x
  \leq |B|^{1/2} \Big( \int_B \E \Big\| \int_{\Gamma (x)} a \D W \Big\|^2 \D x \Big)^{1/2} \leq 1 .
\end{equation*}

Furthermore, for (equivalence classes of) functions $g: \R^{n+1}_+ \to X$ such that
$1_{\Gamma (x;r)}g \in L^2(\R^{n+1}_+) \otimes X$ for every $r > 0$ and almost every $x\in\R^n$ we define
\begin{equation*}
  \| g \|_{T^{\infty}(X)} = \sup_B \Big( \fint_B \E \Big\| \int_{\Gamma (x;r_B)}
  g \D W \Big\| ^2 \D x \Big)^{1/2} < \infty ,
\end{equation*}
where the supremum is taken over all balls $B\subset\R^n$. 
\begin{definition*}
The tent space $T^{\infty}(X)$ is defined as the completion under
$\| \cdot \|_{T^{\infty}(X)}$ of the space of (equivalence classes of) functions 
$g: \R^{n+1}_+ \to X$ such that
$1_{\Gamma (x;r)}g \in L^2(\R^{n+1}_+) \otimes X$ for every $r > 0$ and almost every $x\in\R^n$
and for which $\| g \|_{T^\infty (X)} < \infty$.
\end{definition*}

\subsection*{The atomic decomposition.}

In an atomic decomposition, we aim to express
a $T^1(X)$ function as an infinite sum of (multiples of) atoms. 
The original proof for scalar-valued tent spaces
by Coifman, Meyer and Stein \cite[Theorem 1 (c)]{CMSTENTSPACES} rests on a lemma that allows one to
exchange integration in the upper half-space with ``double integration'', which is something
unthinkable when ``double integration'' consists of both standard and stochastic integration.
The following argument provides a more geometrical reasoning. We start with a covering lemma:

\begin{lem}
\label{balls}
  Suppose that an open set $E\subset\R^n$ has finite measure. Then there exist disjoint balls $B^j\subset E$
  such that
  \begin{equation*}
    \widehat{E} \subset \bigcup_{j\geq 1} \widehat{5B^j} .
  \end{equation*} 
  \begin{proof}
    We start by writing $d_1 = \sup_{B\subset E} r_B$ and choosing a ball $B^1\subset E$ with radius
    $r_1 > d_1 / 2$. Then we proceed inductively: Suppose that balls $B^1,\ldots , B^k$ have been chosen and
    write
    \begin{equation*}
      d_{k+1} = \sup \{ r_B : B\subset E , B\cap B^j = \emptyset , j=1,\ldots ,k \} .
    \end{equation*}
    If possible, we choose $B^{k+1}\subset E$ with radius $r_{k+1} > d_{k+1} / 2$ so that
    $B^{k+1}$ is disjoint from all $B^1, \ldots , B^k$.
    Let then $(y,t)\in \widehat{E}$. In order to show that $B(y,t) \subset 5B^j$
    for some $j$ we note that $B(y,t)$ has to intersect some $B^j$: Indeed, if
    there are only finitely many balls $B^j$, then $y\in\overline{B^j}$ for some $j$. On the other hand,
    if there are infinitely many balls $B^j$ and they are all disjoint from $B(y,t)$, then
    $r_j > d_j / 2 > t / 2$ and $E$ has infinite measure, which is a contradiction. Thus there exists
    a $j$ for which $B(y,t) \cap B^j \neq \emptyset$ and so $B(y,t) \subset 5B^j$ because
    $t\leq d_j \leq 2r_j$ by construction. 
  \end{proof}
\end{lem}

Given a $0 < \lambda < 1$, we define the extension of a measurable set $E\subset\R^n$ by
\begin{equation*}
  E_\lambda^* = \{ x\in\R^n : M1_E(x) > \lambda \} .
\end{equation*}
Here $M$ is the Hardy--Littlewood maximal operator assigning the maximal function
\begin{equation*}
  Mf(x) = \sup_{B\ni x} \fint_B |f(y)| \D y , \quad x\in\R^n ,
\end{equation*}
to every locally integrable real-valued $f$. Note that the lower 
semicontinuity of $Mf$ guarantees that $E_\lambda^*$ is open while
the weak $(1,1)$ inequality for the maximal operator 
assures us that $|E_\lambda^*| \leq \lambda^{-1} |E|$.

We continue by constructing sectors opening in finite number of directions of our choice. To do this, we
fix vectors $v_1,\ldots , v_N$ in the unit sphere $\SP^{n-1}$ of $\R^n$ such that
\begin{equation*}
  \max_{1\leq m \leq N} v \cdot v_m \geq \frac{\sqrt{3}}{2}
\end{equation*}
for every $v\in \SP^{n-1}$. In other words, every $v\in \SP^{n-1}$ makes an angle of no more than
$30^\circ$ with one of $v_m$'s. We write 
\begin{equation*}
  S_m = \Big\{ v\in \SP^{n-1} : v \cdot v_m \geq \frac{\sqrt{3}}{2} \Big\}
\end{equation*} 
and observe that
the angle between two $v, v' \in S_m$ is at most $60^\circ$, i.e. 
$v \cdot v' \geq \frac{1}{2}$. Consequently, $| v - v'| \leq 1$.

For every $x\in\R^n$ and $t > 0$, write 
\begin{equation*}
  R_m(x,t) = \Big\{ y\in B(x,t) : \frac{y-x}{|y-x|} \in S_m \text{ or } y = x \Big\}
\end{equation*}
for the sector opening from $x$ in the direction of $v_m$. For any two 
$y,y'\in R_m(x,t)$, 
the angle between $y-x$ and $y'-x$ is at most $60^\circ$ (when $y$ and $y'$ are different from $x$),
implying that $|y-y'| \leq t$. Hence the proportion
of $R_m(x,t)$ in $B(y,t)$ for any $y\in R_m(x,t)$ is a dimensional constant, in symbols,
\begin{equation*}
  \frac{|R_m(x,t)|}{|B(y,t)|} = c(n), \quad y\in R_m(x,t) .
\end{equation*}
For every $0 < \lambda < c(n)$ it thus holds that
$M1_{R_m(x,t)} > \lambda$ in $B(y,t)$ whenever $y\in R_m(x,t)$. 
Writing $E^* = E_{c(n)/2}^*$ we have now proven the following:

\begin{lem}
\label{sector}
If $E\subset\R^n$ is measurable and $y\in R_m(x,t)\subset E$, then $B(y,t)\subset E^*$.
\end{lem}

Note that the next lemma follows easily when $n=1$ and holds even without the extension.
Indeed, if $E$ is an open interval in $\R$ and $x\in E$, then one can choose $x_1$ and $x_2$ to be the
endpoints of $E$ and obtain $\Gamma (x) \setminus \widehat{E} \subset \Gamma (x_1) \cup \Gamma (x_2)$.
On the other hand, for $n\geq 2$ the extension is necessary, which can be seen already by taking $E$
to be an open ball.

\begin{lem}
\label{cones}
  Suppose that an open set $E\subset\R^n$ has finite measure. Then for every $x\in E$ there exist 
  $x_1,\ldots , x_N\in\partial E$, with $N$ depending only on the dimension $n$, such that
  \begin{equation*}
    \Gamma (x) \setminus \widehat{E^*} \subset \bigcup_{m=1}^N \Gamma (x_m) .
  \end{equation*}
  \begin{proof}
    For every $1 \leq m \leq N$ we may pick $x_m\in\partial E$ in such a manner that
    \begin{equation*}
      \frac{x_m - x}{|x_m - x|} \in S_m
    \end{equation*}
    and $|x_m - x|$, which we denote by $t_m$, is minimal (while positive, since $E$ is open). 
    In other words, $R_m(x,t_m) \subset E$.
    We need to show that for every $(y,t)\in\Gamma (x) \setminus \widehat{E^*}$ 
    the point $y$ is less than $t$ away from
    one of the $x_m$'s. Thus, let $(y,t)\in\Gamma (x) \setminus \widehat{E^*}$, which translates to
    $|x-y| < t$ and $B(y,t)\not\subset E^*$. 
    
    Consider first the case of $y$ not belonging to any 
    $R_m(x,t_m)$.
    Then for some $m$,
    \begin{equation*}
      \frac{y-x}{|y-x|} \in S_m \quad \text{and} \quad |y-x| \geq t_m .
    \end{equation*}
    Now the point
    \begin{equation*}
      z = t_m \frac{y-x}{|y-x|} + x
    \end{equation*}
    sits in the line segment connecting $x$ and $y$ and satisfies $|z-x| = t_m$. Hence the calculation
    \begin{align*}
      |y-x_m| &\leq |y-z|+|z-x_m| \\ 
      &= |y-z|+t_m\Big| \frac{z-x}{t_m} - \frac{x_m-x}{t_m} \Big| \\
      &=|y-z|+|z-x| \Big| \frac{z-x}{|z-x|} - \frac{x_m-x}{|x_m-x|} \Big| \\
      &\leq |y-z|+|z-x| \\
      &= |y-x| < t ,
    \end{align*}
    where we used the fact that $| v - v'| \leq 1$ for any two $v, v' \in S_m$,
    shows that $(y,t)\in \Gamma (x_m)$.
    
    On the other hand, if $y\in R_m(x,t_m)$ for some $m$, then
    $|y-x_m| \leq t_m$, since the diameter of $R_m(x,t_m)$ does not exceed $t_m$. Also
    $B(y,t_m)\subset E^*$ by Lemma \ref{sector} so that $t_m < t$ since $B(y,t)\not\subset E^*$,
    which shows that $(y,t)\in \Gamma (x_m)$.
  \end{proof}
\end{lem}

We are now ready to state and prove the atomic decomposition for $T^1(X)$ functions.

\begin{thm}
  For every function $f$ in $T^1(X)$ 
  there exist countably many atoms $a_k$ and real numbers $\lambda_k$ such that
  \begin{equation*}
    f = \sum_k \lambda_k a_k \quad \text{and} \quad
    \sum_k |\lambda_k | \lesssim \| f \|_{T^1(X)} .
  \end{equation*}
  \begin{proof}
    Let $f$ be a function in $T^1(X)$ and write
    \begin{equation*}
      E_k = \Big\{ x\in\R^n : \Big( \E \Big\| \int_{\Gamma (x)} f \D W \Big\|^2 \Big)^{1/2} > 2^k \Big\}
    \end{equation*}
    for each integer $k$. By Lemma \ref{measurability}, each $E_k$ is open.
    For each $k$, apply Lemma \ref{balls} to the open set $E^*_k$ in order 
    to get disjoint balls $B_k^j\subset E^*_k$ for which
    \begin{equation*}
      \widehat{E_k^*} \subset \bigcup_{j\geq 1} \widehat{5B_k^j} .
    \end{equation*}
    Further, for each of these covers, take a (rough) partition of unity, that is, a collection of
    functions $\chi_k^j$ for which 
    \begin{equation*}    
    0\leq \chi_k^j \leq 1, \quad 
    \sum_{j=1}^{\infty} \chi_k^j = 1 \text{ on } \widehat{E_k^*} \quad \text{and} \quad 
    \supp \chi_k^j \subset \widehat{5B_k^j}.
    \end{equation*}
    For instance, one can define $\chi_k^1$ as the indicator of $\widehat{5B_k^1}$ and $\chi_k^j$ 
    for $j\geq 2$ as the indicator of 
    \begin{equation*}    
    \widehat{5B_k^j} \setminus \bigcup_{i=1}^{j-1} \widehat{5B_k^i} .
    \end{equation*}
    Write $A_k = \widehat{E_k^*}\setminus\widehat{E_{k+1}^*}$.
    We are now in the position to decompose $f$ as
    \begin{equation*}
      f = \sum_{k\in\Z} 1_{A_k} f
      = \sum_{k\in\Z} \sum_{j\geq 1} \chi_k^j 1_{A_k} f
      = \sum_{k\in\Z} \sum_{j\geq 1} \lambda_k^j a_k^j ,
    \end{equation*}
    where
    \begin{equation*}
      \lambda_k^j = |5B_k^j|^{1/2} \Big( \int_{5B_k^j} \E \Big\| 
      \int_{\Gamma (x) \cap A_k} f \D W \Big\|^2 \D x \Big)^{1/2} .
    \end{equation*}
    Observe, that $a_k^j = \chi_k^j 1_{A_k} f / \lambda_k^j$ is an atom supported in $\widehat{5B_k^j}$.
    
    It remains to estimate the sum of $\lambda_k^j$'s.    
    For $x\not\in E_{k+1}$ we have
    \begin{equation*}
      \E \Big\| \int_{\Gamma (x) \cap A_k} f \D W \Big\|^2 \D x \leq 4^{k+1}
    \end{equation*}
    by the definition of $E_{k+1}$.
    The cones at points $x\in E_{k+1}$ are the problematic ones and so
    in order to estimate $\lambda_k^j$'s, we need to exploit the fact that $1_{A_k}f$ vanishes
    on $\widehat{E_{k+1}^*}$.
    Let $x\in E_{k+1}$ and use Lemma \ref{cones} to pick $x_1,\ldots , x_N \in \partial E_{k+1}$,
    where $N\leq c'(n)$,
    such that
    \begin{equation*}
      \Gamma (x) \setminus \widehat{E_{k+1}^*} \subset \bigcup_{m=1}^N \Gamma (x_m) .
    \end{equation*}
    Now $x_1, \ldots , x_N \not\in E_{k+1}$ which allows us to estimate
    \begin{equation*}
      \E \Big\| 
      \int_{\Gamma (x) \cap A_k} f \D W \Big\|^2
      \leq \Big( \sum_{m=1}^N \Big( \E \Big\| \int_{\Gamma (x_m)} f \D W \Big\|^2 \Big) ^{1/2} \Big)^2
      \leq N^2 4^{k+1} .
    \end{equation*}
    Hence, integrating over $5B_k^j$ we obtain
    \begin{equation*}
      \int_{5B_k^j} \E \Big\| 
      \int_{\Gamma (x) \cap A_k} f \D W \Big\|^2 \D x
      \leq |5B_k^j| c'(n)^2 4^{k+1} .
    \end{equation*}
    Consequently,
    \begin{align*}
      \sum_{k\in\Z} \sum_{j\geq 1} \lambda_k^j &\leq c'(n) \sum_{k\in\Z} 2^{k+1} \sum_{j\geq 1} |5B_k^j| \\
      &\leq c'(n)5^n\sum_{k\in\Z} 2^{k+1}|E_k^*| \\ 
      &\leq c'(n)\lambda (n)^{-1} 5^n \sum_{k\in\Z} 2^{k+1}|E_k| \\
      &\leq c'(n)\lambda (n)^{-1} 5^n \| f \|_{T^1(X)} .
    \end{align*}
  \end{proof}
\end{thm}

It is perhaps surprising that the UMD assumption is not needed for the atomic decomposition.

\subsection*{Embedding $T^1(X)$ into $H^1(\R^n ; \gamma (X))$ and $T^{\infty}(X)$ into $\textup{BMO}(\R^n ; \gamma (X))$.}

Armed with the atomic decomposition we proceed to the embeddings.
Suppose that a smooth function $\psi : [0,\infty ) \to \R$ satisfies 
$1_{[0,1)} \leq |\psi | \leq 1_{[0,\alpha )}$
for some $\alpha > 2$ and has $\int_{\R^n} \psi (|x|) \D x = 0$.
For functions $f : \R^{n+1}_+ \to X$ we define
\begin{equation*}
  J_\psi f(x;y,t) = \psi \Big( \frac{|x-y|}{t} \Big) f(y,t), \quad x\in\R^n , (y,t)\in\R^{n+1}_+ ,
\end{equation*}
and note immediately that $\int_{\R^n} J_\psi f (x) \D x = 0$.

Recall also that functions in the Hardy space $H^1(\R^n ; \gamma (X))$ are composed of
atoms $A: \R^n \to \gamma (X)$ 
each of which is supported on a ball $B\subset\R^n$, has zero integral and satisfies
\begin{equation*}
  \int_B \E \Big\| \int_{\R^{n+1}_+} A(x;y,t) \D W(y,t) \Big\|^2 \D x \leq \frac{1}{|B|} .
\end{equation*}
For further references, see Blasco \cite{BLASCOHARDY} and Hyt\"onen \cite{HYTONENHARDY}.

\begin{thm}
  Suppose that $X$ is UMD. Then $J_\psi$ embeds $T^1(X)$ into $H^1(\R^n ; \gamma (X))$
  and $T^{\infty}(X)$ into $\textup{BMO}(\R^n ; \gamma (X))$.
\end{thm}
\begin{proof}
We argue that $J_\psi$ takes $T^1(X)$ atoms to (multiples of) $H^1(\R^n ; \gamma (X))$ atoms.
If a $T^1(X)$ atom $a$ is supported in $\widehat{B}$ for some ball $B\subset\R^n$, then
$J_\psi a$ is supported in $\alpha B$ and $\int J_\psi a = 0$.
Moreover, since $X$ is UMD, we may use the equivalence of
$T^2(X)$ norms (Theorem \ref{normequiv}) and write
\begin{equation*}
  \int_{\alpha B} \E \Big\| \int_{\R^{n+1}_+} \psi 
  \Big( \frac{|x-y|}{t} \Big) a(y,t) \D W(y,t) \Big\| ^2 \D x
  \lesssim \int_B \E \Big\| \int_{\Gamma (x)} a \D W \Big\| ^2 \D x \leq \frac{1}{|B|} .
\end{equation*}
The boundedness of $J_\psi$ from $T^1(X)$ to $H^1(\R^n ; \gamma (X))$ follows.
In addition, since $1_{[0,1)} \leq |\psi |$, it follows that
$\| f \|_{T^1(X)} \leq \| J_\psi f \|_{L^1(\R^n ; \gamma (X))} \leq \| J_\psi f \|_{H^1(\R^n ; \gamma (X))}$
and so
$J_\psi$ is also bounded from below.

To see that $J_\psi$ maps $T^{\infty}(X)$ boundedly into $\textup{BMO}(\R^n ; \gamma (X))$, we need to show that
\begin{equation*}
  \Big( \fint_B \E \Big\| \int_{\R^{n+1}_+} \Big( J_\psi g(x;y,t) -
  \fint_B J_\psi g(z;y,t) \D z \Big) \D W(y,t) \Big\|^2 \D x \Big)^{1/2} \lesssim \| g \|_{T^{\infty}(X)}
\end{equation*}
for all balls $B\subset \R^n$.
We partition the upper half-space into 
$\R^n \times (0,r_B)$ and the sets
$A_k = \R^n \times [2^{k-1}r_B,2^kr_B)$ for positive integers $k$ 
and study each piece separately.

On $\R^n \times (0,r_B)$ one has
\begin{align*}
  &\Big( \fint_B \E \Big\| \int_{\R^n\times (0,r_B)} \psi \Big( \frac{|z-y|}{t} \Big)
  g(y,t) \D W(y,t) \Big\|^2 \D z \Big)^{1/2} \\
  &\leq \Big( \fint_B \E \Big\| \int_{\Gamma_\alpha (x;r_B)} g \D W 
  \Big\|^2 \D x \Big)^{1/2} 
  \lesssim \| g \|_{T^{\infty}}
\end{align*}
since $|\psi | \leq 1_{[0,\alpha )}$ and the $T^2(X)$ norms are comparable (Theorem \ref{normequiv}).
Furthermore, as one can justify by approximating $\psi$ with simple functions, we have
\begin{align*}
  &\Big( \E \Big\| \int_{\R^n\times (0,r_B)} g(y,t) 
  \fint_B\psi \Big( \frac{|z-y|}{t} \Big) \D z
  \D W(y,t) \Big\|^2 \Big)^{1/2} \\
  &\leq \Big( \fint_B \E \Big\| \int_{\R^n\times (0,r_B)} \psi \Big( \frac{|z-y|}{t} \Big)
  g(y,t) \D W(y,t) \Big\|^2 \D z \Big)^{1/2} ,
\end{align*}
which can be estimated from above by $\| g \|_{T^{\infty}}$, as above.

For each $k$ and $x \in B$, we claim that
\begin{equation*}
  \Big| \fint_B \Big( \psi \Big( \frac{|x-y|}{t} \Big) - \psi \Big( \frac{|z-y|}{t} \Big) \Big)
  \D z \Big| \lesssim 2^{-k} 1_{\Gamma_{\alpha + 2}(x)}(y,t) ,
\end{equation*}
whenever $(y,t) \in A_k$.
Indeed, if $(y,t) \in A_k \cap \Gamma_{\alpha + 2}(x)$, we may use the fact that
\begin{equation*}
  \Big| \psi \Big( \frac{|x-y|}{t} \Big) - \psi \Big( \frac{|z-y|}{t} \Big) \Big| \lesssim
  \sup |\psi'| \frac{|x-z|}{t} \lesssim \frac{r_B}{2^k r_B} = 2^{-k}
\end{equation*}
for all $z \in B$,
while for $(y,t)\in A_k \setminus \Gamma_{\alpha + 2}(x)$ we have 
$|y-x| \geq (\alpha + 2)t \geq \alpha t + 2r_B$ so that $|y-z| \geq |y-x| - |x-z| \geq \alpha t$
for each $z\in B$, which results in
\begin{equation*}
  \int_B \Big( \psi \Big( \frac{|x-y|}{t} \Big) - \psi \Big( \frac{|z-y|}{t} \Big) \Big)
  \D z = 0 . 
\end{equation*}

This gives
\begin{align*}
  &\Big( \fint_B \E \Big\| \int_{A_k}
  \frac{g(y,t)}{|B|} \int_B \Big( \psi \Big( \frac{|x-y|}{t} \Big) - \psi \Big( \frac{|z-y|}{t} \Big) \Big)
  \D z \D W(y,t) \Big\|^2 \D x \Big)^{1/2} \\
  &\leq 2^{-k} \Big( \fint_B \E \Big\| \int_{A_k \cap \Gamma_{\alpha + 2} (x)}
  g \D W \Big\| ^2 \D x \Big)^{1/2} .
\end{align*}
But every $A_k \cap \Gamma_{\alpha + 2} (x)$ with $x\in B$ is contained in any $\Gamma_{\alpha + 6} (z)$ with
$z\in 2^kB$. Indeed, for all $(y,t) \in A_k \cap \Gamma_{\alpha + 2} (x)$ we have
\begin{equation*}
  |y-z| \leq |y-x| + |x-z| \leq (\alpha + 2) t + (2^k + 1)r_B \leq (\alpha + 6)t .
\end{equation*}
Hence
\begin{equation*}
  \fint_B \E \Big\| \int_{A_k \cap \Gamma_{\alpha + 2} (x)}
  g \D W \Big\| ^2 \D x
  \leq \fint_{2^kB} \E \Big\| \int_{\Gamma_{\alpha + 6} (z)}
  g \D W \Big\| ^2 \D z .
\end{equation*}
Summing up, we obtain
\begin{align*}
  &\sum_{k=1}^{\infty} \Big( \fint_B \E \Big\| \int_{A_k}
  g(y,t) \fint_B \Big( \psi \Big( \frac{|x-y|}{t} \Big) - \psi \Big( \frac{|z-y|}{t} \Big) \Big)
  \D z \D W(y,t) \Big\|^2 \D x \Big)^{1/2} \\
  &\leq \sum_{k=1}^{\infty} 2^{-k} \Big( 
  \fint_{2^kB} \E \Big\| \int_{\Gamma_{\alpha + 6} (z)}
  g \D W \Big\| ^2 \D z \Big)^{1/2} 
  \lesssim \| g \|_{T^{\infty}(X)} .
\end{align*}

To see that $\| g \|_{T^{\infty}(X)} \lesssim \| J_\psi g \|_{\textup{BMO}(\R^n ; \gamma (X))}$ it suffices to
fix a ball $B\subset\R^n$ and show,
that for every $x\in B$ we have
\begin{equation*}
  1_{\Gamma (x;r_B)} (y,t) \leq \Big| \psi \Big( \frac{|x-y|}{t} \Big) 
  - \fint_{(\alpha +2)B}
  \psi \Big( \frac{|z-y|}{t} \Big) \D z \Big| ,
\end{equation*}
since this gives us
\begin{align*}
  &\fint_B \E \Big\| \int_{\Gamma (x;r_B)} g \D W \Big\|^2 \D x \\
  &\leq \fint_B \E \Big\| \int_{\R^{n+1}_+} g(y,t) 
  \Big( \psi \Big( \frac{|x-y|}{t} \Big) - 
  \fint_{(\alpha + 2)B} \psi \Big( \frac{|z-y|}{t} \Big) \D z \Big) \Big\|^2 \D x \\
  &\leq (\alpha + 2)^n \| J_\psi g \|_{\textup{BMO}(\R^n ; \gamma (X))} .
\end{align*}
Now that $1_{[0,1)} \leq | \psi |$ and $\int_{\R^n} \psi (|x|) \D x = 0$, it is enough to prove for
a fixed $x\in B$, that 
\begin{equation*}
  \supp \psi \Big( \frac{|\cdot - y|}{t} \Big) \subset (\alpha + 2)B
\end{equation*}
for every 
$(y,t)\in\Gamma (x;r_B)$, i.e. that $B(y,\alpha t) \subset (\alpha + 2)B$ whenever 
$|x-y| < t < r_B$. This is indeed true, as every $z\in B(y,\alpha t)$ satisfies
\begin{equation*}
  |z - x| \leq |z - y| + |y - x| < (\alpha + 1)r_B .
\end{equation*}
We have established that, also in this case, $J_\psi$ is bounded from below. 
\end{proof}

It follows that different $T^1(X)$ norms are equivalent in the sense that whenever
$1_{[0,1)} \leq |\phi | \leq 1_{[0,\alpha )}$ for some $\alpha > 1$, we can take smooth
$\psi : [0,\infty ) \to \R$ with $|\phi | \leq |\psi | \leq 1_{[0,2\alpha )}$ to obtain
\begin{equation*}
  \| f \|_{T^1(X)} \leq \| J_\phi f \|_{L^1(\R^n ; \gamma (X))}
  \leq \| J_\psi f \|_{L^1(\R^n ; \gamma (X))} \leq \| J_\psi f \|_{H^1(\R^n ; \gamma (X))}
  \lesssim \| f \|_{T^1(X)} .
\end{equation*}

To identify $T^1(X)$ as a complemented subspace of $H^1(\R^n ; \gamma (X))$ we
define a projection first on the level of test functions.
Let us write 
\begin{equation*}
  T(X) = \{ f:\R^{n+1}_+ \to X : \, 1_{\Gamma (x)}f \in L^2(\R^{n+1}_+) \otimes X
  \text{ for almost every } x\in\R^n \}
\end{equation*}
and
\begin{align*}
  S(\gamma (X)) = \textup{span}\, \{ &F:\R^n \times \R^{n+1}_+ \to X : \,
  F(x;y,t) = \Psi (x;y,t)f(y,t) \\ &\text{for some } \Psi\in L^\infty (\R^n \times \R^{n+1}_+)
  \text{ and } f\in T(X) \} .
\end{align*} 
Observe, that $J_\psi$ maps $T(X)$ into $S(\gamma (X))$
and that $S(\gamma (X))$ intersects $L^p(\R^n ; \gamma (X))$ densely for all $1 < p < \infty$ and likewise for
$H^1(\R^n ; \gamma (X))$.

For $F$ in $S(\gamma (X))$ we define
\begin{equation*}
  (N_\psi F)(x;y,t) = \psi \Big( \frac{|x-y|}{t} \Big) \frac{1}{c_\psi t^n}
  \int_{\R^n} \psi \Big( \frac{|z-y|}{t} \Big) F(z;y,t) \D z ,
\end{equation*}
where $c_\psi = \int_{\R^n} \psi (|x|)^2 \D x$. Now $N_\psi$ is a projection and
satisfies $N_\psi J_\psi = J_\psi$. Also,
for every $F\in S(\gamma (X))$ we find an $f\in T(X)$ so that
$N_\psi F = J_\psi f$, namely
\begin{equation*}
  f(y,t) = \frac{1}{c_\psi t^n} 
  \int_{\R^n} \psi \Big( \frac{|z-y|}{t} \Big) F(z;y,t) \D z , \quad (y,t)\in\R^{n+1}_+ .
\end{equation*}

\begin{thm}
\label{endpointembedding}
  Suppose that $X$ is UMD. Then $N_\psi$ extends to a bounded projection on $H^1(\R^n ; \gamma (X))$
  and $J_\psi$ extends to an isomorphism from $T^1(X)$ onto the image of $H^1(\R^n ; \gamma (X))$ 
  under $N_\psi$.
  \begin{proof}
    Let $1 < p < \infty$.
    For simple $L^2(\R^{n+1}_+) \otimes X$ -valued functions $F$ defined on $\R^n$ the mapping
    $(y,t) \mapsto F(\cdot ; y,t) : \R^{n+1}_+ \to L^p(\R^n ; X)$ is
    in $L^2(\R^{n+1}_+) \otimes L^p(\R^n ; X)$ and
    we may express $N_\psi$ using the averaging operators as
    \begin{equation*}
      (N_\psi F)( \cdot ; y,t) = A_{y,t}^{\psi} (F(\cdot ; y,t)) .
    \end{equation*}
    Since $X$ is UMD, Stein's inequality guarantees $\gamma$-boundedness for the range of the
    strongly $L^p(\R^n ; X)$-measurable function 
    $(y,t)\mapsto A_{y,t}^\psi$, and so by Lemma \ref{gammabounded},
    \begin{equation*}    
      \E \Big\| \int_{\R^{n+1}_+} A_{y,t}^\psi (F(\cdot ; y,t)) \D W(y,t) \Big\|_{L^p(\R^n ; X)}^p
      \lesssim \E \Big\| \int_{\R^{n+1}_+} F(\cdot ; y,t) \D W(y,t) \Big\|_{L^p(\R^n ; X)}^p .
    \end{equation*}
    In other words, $\| N_\psi F \|_{L^p(\R^n ; \gamma (X))}^p \lesssim \| F \|_{L^p(\R^n ; \gamma (X))}^p$.
    
    We wish to define a suitable $\mathcal{L}(\gamma (X))$-valued kernel $K$ that allows us to 
    express $N_\psi$ as a Calder\'on--Zygmund operator
    \begin{equation*}
      N_\psi F(x) = \int_{\R^n} K(x,z)F(z)\D z, \quad F\in L^p(\R^n ; \gamma (X)) .
    \end{equation*}
    For distinct $x,z\in\R^n$ and we define $K(x,z)$ as multiplication by
    \begin{equation*}
      (y,t) \mapsto \psi \Big( \frac{|x-y|}{t} \Big) \frac{1}{c_\psi t^n} 
      \psi \Big( \frac{|z-y|}{t} \Big) ,
    \end{equation*}
    and so
    \begin{equation*}
      \| K(x,z) \|_{\mathcal{L}(\gamma (X))} 
      = \sup_{(y,t)\in\R^{n+1}} \Big| \psi \Big( \frac{|x-y|}{t} \Big)
      \frac{1}{c_\psi t^n} \psi \Big( \frac{|z-y|}{t} \Big) \Big| .
    \end{equation*}
    For $|x-z| > \alpha t$ we have
    \begin{equation*}
      \psi \Big( \frac{|x-y|}{t} \Big)
      \frac{1}{c_\psi t^n} \psi \Big( \frac{|z-y|}{t} \Big) = 0 
    \end{equation*}
    while $|x-z| \leq \alpha t$ guarantees that
    \begin{equation*}
      \Big| \psi \Big( \frac{|x-y|}{t} \Big)
      \frac{1}{c_\psi t^n} \psi \Big( \frac{|z-y|}{t} \Big) \Big| \leq \frac{1}{c_\psi t^n} 
      \leq \frac{\alpha ^n}{c_\psi |x-z|^n} . 
    \end{equation*}
    Hence
    \begin{equation*}
      \| K(x,z) \|_{\mathcal{L}(\gamma (X))} \lesssim \frac{1}{|x-z|^n} .
    \end{equation*}
    Similarly,
    \begin{equation*}
      \| \nabla _x K(x,z) \|_{\mathcal{L}(\gamma (X))}
      = \sup_{(y,t)\in\R^{n+1}_+} \Big| \psi ' \Big( \frac{|x-y|}{t} \Big) \frac{1}{c_\psi t^{n+1}}
      \psi \Big( \frac{|z-y|}{t} \Big) \Big| \lesssim \frac{1}{|x-z|^{n+1}} .
    \end{equation*}
    Thus $K$ is indeed a Calder\'on--Zygmund kernel.
    
    Now $\int_{\R^n} \psi (|x|) \D x = 0$ implies that 
    $\int_{\R^n} N_\psi F (x) \D x = 0$ for $F\in H^1(\R^n ; \gamma (X))$, which
    guarantees that $N_\psi$ maps $H^1(\R^n ; \gamma (X))$ into itself (see Meyer and Coifman
    \cite[Chapter 7, Section 4]{COIFMANMEYER}).
  \end{proof}
\end{thm}

We proceed to the question of duality of $T^1(X)$ and $T^\infty (X^*)$.
Assuming that $X$ is UMD, it is both reflexive and K-convex so that the duality
\begin{equation*}
  H^1(\R^n ; \gamma (X))^* \simeq \textup{BMO}(\R^n ; \gamma (X)^*) \simeq \textup{BMO}(\R^n ; \gamma (X^*))
\end{equation*}
holds (recall the discussion in Section 2) 
and we may define the adjoint of $N_\psi$ by $\la F , N_\psi^* G \ra = \la N_\psi F , G \ra$, where
$F\in H^1(\R^n ; \gamma (X))$ and $G\in \textup{BMO}(\R^n ; \gamma (X^*))$. Moreover, as $T^1(X)$ is
isomorphic to the image of $H^1(\R^n ; \gamma (X))$ under $N_\psi$, its dual $T^1(X)^*$ is isomorphic
to the image of $\textup{BMO}(\R^n ; \gamma (X^*))$ under the adjoint $N_\psi^*$ and the question arises whether
the latter is isomorphic to $T^\infty (X^*)$. For $J_\psi$ to give this isomorphism 
(and to be onto) one could try and
follow the proof strategy of the case $1 < p < \infty$ and give an explicit definition of $N_\psi^*$
on a dense subspace of $\textup{BMO}(\R^n ; \gamma (X^*))$. Even though 
the properties of the kernel $K$ of $N_\psi$ guarantee that
$N_\psi^*$ formally agrees with $N_\psi$ on $L^p(\R^n ; \gamma (X^*))$, it is problematic to
find suitable dense subspaces of $\textup{BMO}(\R^n ; \gamma (X^*))$.

In order to address these issues in more detail, we specify another pair of test function classes,
namely
\begin{align*}
  \widetilde{T}(X) = \{ g:\R^{n+1}_+ \to X : \, &1_{\Gamma (x;r)}g \in L^2(\R^{n+1}_+) \otimes X
  \text{ for every } r > 0 \\ &\text{and for almost every } x\in\R^n \}
\end{align*}
and
\begin{align*}
  \widetilde{S}(\gamma (X)) = \textup{span}\, \{ &G:\R^n \times \R^{n+1}_+ \to X : \,
  G(x;y,t) = \Psi (x;y,t)g(y,t) \text{ for some} \\
  &\Psi\in L^\infty (\R^n \times \R^{n+1}_+)
  \text{ and } g\in \widetilde{T}(X) \} \, / \, \{ \text{constant functions} \} .
\end{align*}

Since $\int_{\R^n} \psi (|x|) \D x = 0$, the projection $N_\psi$ is well-defined on $\widetilde{S}(\gamma (X))$.
Moreover, given any $G\in \widetilde{S}(\gamma (X))$ we can write
\begin{equation*}
  g(y,t) = \frac{1}{c_\psi t^n} \int_{\R^n} \psi \Big( \frac{|z-y|}{t} \Big) G(z;y,t) \D z 
\end{equation*}
to define a function $g\in \widetilde{T}(X)$ for which $N_\psi G = J_\psi g$. But 
$\widetilde{S}(\gamma (X))$ has only weak*-dense
intersection with $\textup{BMO}(\R^n ; \gamma (X))$ (recall that $X \simeq X^{**}$). 
Nevertheless,
$J_\psi$ is an isomorphism from $T^{\infty}(X)$ onto the closure of the image of 
$\widetilde{S}(\gamma (X)) \cap \textup{BMO}(\R^n ; \gamma (X))$ under $N_\psi$. 
It is not clear whether test functions are dense in the closure of their image
under the projection.

The following relaxed duality result is still valid:

\begin{thm}
  Suppose that $X$ is UMD. 
  Then $T^{\infty}(X^*)$ isomorphic to a norming subspace of $T^1(X)^*$ and its action is realized
  for functions $f\in T^1(X)$ and $g\in T^{\infty}(X^*)$ via 
  \begin{equation*}
    \la f,g \ra = c \int_{\R^{n+1}_+} \la f(y,t),g(y,t) \ra \frac{\D y \D t}{t} ,
  \end{equation*}
  where $c$ depends on the dimension $n$.
  \begin{proof}
    Fix a smooth $\psi : [0,\infty ) \to \R$ such that $1_{[0,1)} \leq |\psi | \leq 1_{[0,\alpha )}$
    for some \mbox{$\alpha > 2$} and $\int_{\R^n} \psi (|x|) \D x = 0$.
    By Theorem \ref{endpointembedding}, $T^1(X)$ is isomorphic to the image of $H^1(\R^n ; \gamma (X))$ under
    $N_\psi$, from which it follows that the dual $T^1(X)^*$ is isomorphic to
    the image of $\textup{BMO}(\R^n ; \gamma (X^*))$ under the adjoint projection $N_\psi^*$, 
    which formally agrees with $N_\psi$. 
    The space $T^{\infty}(X^*)$, on the other hand, is 
    isomorphic to the closure of the image of 
    $\widetilde{S}(\gamma (X^*)) \cap \textup{BMO}(\R^n ; \gamma (X^*))$ under $N_\psi$ in
    $\textup{BMO}(\R^n ; \gamma (X^*))$ and hence is a closed subspace of $T^1(X)^*$.
    We can pair a function $f\in T^1(X)$ with a function $g\in T^{\infty}(X^*)$ using the pairing
    of $J_\psi f$ and $J_\psi g$ and
    the atomic decomposition of $f$ to get:
    \begin{align*}
      \la f,g \ra = \sum_k \la J_\psi a_k , J_\psi g \ra 
      &= \sum_k \lambda_k \int_{\R^n} \int_{\R^{n+1}_+} \psi \Big( \frac{|x-y|}{t} \Big) ^2
      \la a_k(y,t),g(y,t) \ra \frac{\D y \D t}{t^{n+1}} \\      
      &= c_nc_\psi \sum_k \lambda_k
      \int_{\R^{n+1}_+} \la a_k(y,t) , g(y,t) \ra \frac{\D y \D t}{t} \\
      &= c_nc_\psi \int_{\R^{n+1}_+} \la f(y,t) , g(y,t) \ra \frac{\D y \D t}{t} ,
    \end{align*}
    where $c_n$ denotes the volume of the unit ball in $\R^n$. 
    The subspace $L^\infty (\R^n) \otimes L^2(\R^{n+1}_+) \otimes X^*$ 
    is weak*-dense in $\textup{BMO}(\R^n ; \gamma (X^*))$ and hence a norming subspace for
    $H^1(\R^n ; \gamma (X))$.
    As it is contained in $\widetilde{S}(\gamma (X^*)) \cap \textup{BMO}(\R^n ; \gamma (X^*))$, we obtain
    \begin{align*}
      \| f \|_{T^1(X)} &\eqsim \| J_\psi f \|_{H^1(\R^n ; \gamma (X))}
      = \sup_{G} | \la J_\psi f , G \ra | = \sup_{G} | \la N_\psi J_\psi f , G \ra | \\
      &= \sup_{G} | \la J_\psi f , N_\psi^* G \ra |
      \eqsim \sup_{g} | \la J_\psi f , J_\psi g \ra | = \sup_g | \la f , g \ra | ,
    \end{align*}
    where the suprema are taken over 
    $G\in\widetilde{S}(\gamma (X^*)) \cap \textup{BMO}(\R^n ; \gamma (X^*))$ with \linebreak
    $\| G \|_{\textup{BMO}(\R^n ; \gamma (X^*))} \leq 1$ and $g\in T^{\infty}(X^*)$ with 
    $\| g \|_{T^{\infty}(X^*)} \leq 1$.
  \end{proof}
\end{thm}


\bibliographystyle{plain}
\bibliography{JAMStentspaces}

\begin{thebibliography}{10}

\bibitem{BLASCOHARDY}
O.~Blasco.
\newblock Hardy spaces of vector-valued functions: duality.
\newblock {\em Trans. Amer. Math. Soc.}, 308(2):495--507, 1988.

\bibitem{BOURGAINSTEININEQ}
J.~Bourgain.
\newblock Vector-valued singular integrals and the {$H^1$}-{BMO} duality.
\newblock In {\em Probability theory and harmonic analysis ({C}leveland,
  {O}hio, 1983)}, volume~98 of {\em Monogr. Textbooks Pure Appl. Math.}, pages
  1--19. Dekker, New York, 1986.

\bibitem{CLEMENT}
P.~Cl{\'e}ment, B.~de~Pagter, F.~A. Sukochev, and H.~Witvliet.
\newblock Schauder decomposition and multiplier theorems.
\newblock {\em Studia Math.}, 138(2):135--163, 2000.

\bibitem{CMSTENTSPACES}
R.~R. Coifman, Y.~Meyer, and E.~M. Stein.
\newblock Some new function spaces and their applications to harmonic analysis.
\newblock {\em J. Funct. Anal.}, 62(2):304--335, 1985.

\bibitem{HARBOURE}
E.~Harboure, J.~L. Torrea, and B.~E. Viviani.
\newblock A vector-valued approach to tent spaces.
\newblock {\em J. Analyse Math.}, 56:125--140, 1991.

\bibitem{HYTONENHARDY}
T.~Hyt{\"o}nen.
\newblock Vector-valued wavelets and the {H}ardy space {$H^1(\Bbb R^n,X)$}.
\newblock {\em Studia Math.}, 172(2):125--147, 2006.

\bibitem{HVNPCONICAL}
T.~Hyt{\"o}nen, J.~M. A. M.~van Neerven, and P.~Portal.
\newblock Conical square function estimates in {UMD} {B}anach spaces and
  applications to {$H^\infty$}-functional calculi.
\newblock {\em J. Anal. Math.}, 106:317--351, 2008.

\bibitem{HYTONENWEISPARAPRODUCTS}
T.~Hyt{\"o}nen and L.~Weis.
\newblock The {B}anach space-valued {BMO}, {C}arleson’s condition, and
  {P}araproducts.
\newblock {\em Journal of Fourier Analysis and Applications}, 16:495--513.

\bibitem{JANSON}
S.~Janson.
\newblock {\em Gaussian {H}ilbert spaces}, volume 129 of {\em Cambridge Tracts
  in Mathematics}.
\newblock Cambridge University Press, Cambridge, 1997.

\bibitem{WEISUNPUB}
N.~Kalton and L.~Weis.
\newblock The {$H^\infty$}-functional calculus and square function estimates.
\newblock Manuscript in preparation.

\bibitem{TAOMEI}
T.~Mei.
\newblock B{MO} is the intersection of two translates of dyadic {BMO}.
\newblock {\em C. R. Math. Acad. Sci. Paris}, 336(12):1003--1006, 2003.

\bibitem{COIFMANMEYER}
Y.~Meyer and R.~R. Coifman.
\newblock {\em Wavelets}, volume~48 of {\em Cambridge Studies in Advanced
  Mathematics}.
\newblock Cambridge University Press, Cambridge, 1997.
\newblock Calder{\'o}n-Zygmund and multilinear operators, Translated from the
  1990 and 1991 French originals by David Salinger.

\bibitem{GAMMARAD}
J.~M. A. M.~van Neerven.
\newblock {$\gamma$}-radonifying operators---a survey.
\newblock In {\em The {AMSI}-{ANU} {W}orkshop on {S}pectral {T}heory and
  {H}armonic {A}nalysis}, volume~44 of {\em Proc. Centre Math. Appl. Austral.
  Nat. Univ.}, pages 1--61. Austral. Nat. Univ., Canberra, 2010.

\bibitem{JVNSTOCHINT}
J.~M. A. M.~van Neerven and L.~Weis.
\newblock Stochastic integration of functions with values in a {B}anach space.
\newblock {\em Studia Math.}, 166(2):131--170, 2005.

\bibitem{ROSINSKI}
J.~Rosi{\'n}ski and Z.~Suchanecki.
\newblock On the space of vector-valued functions integrable with respect to
  the white noise.
\newblock {\em Colloq. Math.}, 43(1):183--201 (1981), 1980.

\end{thebibliography}

\end{document}